\documentclass[a4paper,reqno,12pt]{amsart}
\usepackage{amssymb}
\usepackage{yhmath}       
\usepackage{bbm}          
\usepackage{mathrsfs}     

\usepackage[T1]{fontenc}
\usepackage{textcomp}
\usepackage{amsmath,amsthm,amssymb,amscd}
\usepackage{type1cm}
\usepackage{txfonts}

\numberwithin{equation}{section}	

\theoremstyle{plain}
\newtheorem{thm}{Theorem}[section]      
\newtheorem*{thm*}{Theorem}

\newtheorem{lem}[thm]{Lemma}

\newtheorem{prop}[thm]{Proposition}

\newtheorem{cor}[thm]{Corollary}

\theoremstyle{definition}

\theoremstyle{remark}
\newtheorem{remark}[thm]{Remark}

\theoremstyle{remark}

\def\geqsl{\geqslant}
\def\leqsl{\leqslant}
\newcommand{\floor}[1]{\lfloor{#1}\rfloor}
\newcommand{\tild}[1]{\widetilde{#1}}

\newcommand{\normalordering}[1]{%
  \ensuremath{{}_{\circ}^{\circ} }\,#1\,{ _{\circ}^{\circ} }}

\renewcommand{\setminus}{\smallsetminus}  

\newcommand{\sgn}[1]{\ensuremath{\operatorname{sgn}{\left( #1 \right)}}}
\newcommand{\tp}[1]{\ensuremath{\operatorname{}^{t}\hspace{-3truept}{#1}}}
\newcommand{\ttp}[1]{\ensuremath{#1^{tt}}}
\newcommand{\Mat}[1]{\ensuremath{\operatorname{Mat}_{#1}}}
\newcommand{\Alt}{\ensuremath{\operatorname{Alt}}}   
   
\newcommand{\bs}[1]{\ensuremath{\boldsymbol{#1}}}
\newcommand{\Pf}[1]{\ensuremath{\operatorname{Pf}\left({#1}\right)}}

\newcommand{\Ad}{\ensuremath{\operatorname{Ad}}}  

\def\c{\gamma}

\def\A{\mathscr A}	%

\renewcommand{\Omega}{\ensuremath{\varOmega}}
\renewcommand{\Xi}{\ensuremath{\varXi}}
\renewcommand{\Theta}{\ensuremath{\varTheta}}

\newcommand{\N}{\ensuremath{\mathbb N}}

\newcommand{\C}{\ensuremath{\mathbb C}}

\newcommand{\gl}{\ensuremath{\mathfrak {gl}}}

\newcommand{\GL}[1]{\ensuremath{\mathrm{GL}_{#1}}}

\newcommand{\SO}{\ensuremath{\mathrm{SO}}}   
   
\renewcommand{\d}{\ensuremath{\operatorname{d}\!}}
\newcommand{\pd}{\ensuremath{\partial}}

\def\<{\langle}
\def\>{\rangle}

\def\c2vec#1#2{ %
   \left[ \begin{smallmatrix} %
           #1 \\ #2  \end{smallmatrix} %
   \right]}

\numberwithin{equation}{section}

\renewcommand{\tp}[1]{\ensuremath{#1^{t}}}
\renewcommand{\Pf}[1]{\ensuremath{\operatorname{Pf}({#1})}}

\newcommand{\PD}{\ensuremath{\mathscr{PD}}}
\newcommand{\K}{\ensuremath{\mathbb{K}}}

\newcommand{\dfac}[2]{\ensuremath{{#1}^{\,\underline{#2}}}}

\title%
[Generating function for invariant differential operators]
{Generating function for $\GL{n}$-invariant \\
 differential operators in the skew Capelli identity}
\author{Takashi Hashimoto}
\address{%
 Department of Informatics, 
 Graduate School of Engineering,
 Tottori University,
 4-101, Koyama-Minami, Tottori, 680-8552, Japan
        }
\email{thashi@ike.tottori-u.ac.jp}
\date{\today}

\begin{document}

\begin{abstract}
Let $\Alt_n$ be the vector space of
all alternating $n \times n$ complex matrices,
on which 
the complex general linear group $\GL{n}$
acts by $x \mapsto g x \tp{g}$. 
The aim of this paper is 
to show that Pfaffian of a certain matrix 
whose entries are 
multiplication operators
or derivations acting on polynomials on $\Alt_n$
provides a generating function for 
the $\GL{n}$-invariant differential operators 
that play a role 
in the skew Capelli identity,
with coefficients the Hermite polynomials. 

\end{abstract}

\maketitle

\section{Introduction}

Let $\Alt_n$ be the vector space consisting 
of all alternating $n \times n$ complex matrices,
and 
$\C[\Alt_n]$ the $\C$-algebra of all polynomial functions on $\Alt_n$.
Then the complex general linear group $\GL{n}$ 
acts on $\Alt_n$ by
\begin{equation}
\label{e:GL_action}
  g.x:=g x \tp{g}  \quad (g \in \GL{n}, x \in \Alt_n),
\end{equation}
from which 
one can define a representation $\pi$ of $\GL{n}$ 
on $\C[\Alt_n]$ by
\begin{equation}
\label{e:GL_rep_on_polynom}
  \pi(g)f(x):=f(g^{-1}.x)   
  \quad  (g \in \GL{n}, f \in \C[\Alt_n]).
\end{equation}
For $x=(x_{i,j})_{i,j=1,\dots,n} \in \Alt_n$,
with $x_{j,i}=-x_{i,j}$,
let $M:=(x_{i,j})_{i,j}$ and $D:=({\pd}_{i,j})_{i,j}$
be the alternating $n \times n$ matrices 
whose $(i,j)$-th entries are given by
the multiplication operator $x_{i,j}$
and the derivation $\pd_{i,j}:=\pd/\pd x_{i,j}$,
respectively.
Then the representation $\d\pi$ of $\gl_n$,  
the Lie algebra of $\GL{n}$,
induced from $\pi$ is given by
\begin{equation}
\label{e:gl_rep_on_polynom}
\d\pi(E_{i,j}) = -\sum_{k=1}^{n} x_{k,j} \pd_{k,i}
  \quad (i,j=1,2,\dots,n)
\end{equation}
where $E_{i,j}$ denotes the matrix unit of size $n \times n$
which is a basis for $\gl_n$.

Let us denote 
the ring of differential operators on $\Alt_n$ 
with polynomial coefficients
by $\PD(\Alt_n)$,
and its subring consisting of 
$\GL{n}$-invariant differential operators
by $\PD(\Alt_n)^{\GL{n}}$.
Moreover,
for a positive integer $n$,
$[n]$ denotes the set $\{ 1,2,\dots,n \}$,
and for a real number $x$,
$\floor{x}$ the greatest integer not exceeding $x$.
Then the following fact is known:
\begin{thm*}
[\cite{HU91}]
 For $k=0,1,\dots,\floor{{n}/{2}}$,
 let
 \begin{equation}
  \varGamma_k:=\sum_{ I \in \binom{[n]}{2k} }
                        \Pf{x_I} \Pf{\pd_I},
     \label{e:skew_capelli}
 \end{equation}
 where the summation is taken over all 
 $I \subset [n]$ such that
 its cardinality is $2k$,
 and $x_{I}$ and $\pd_I$ denote submatrices of $M$ and $D$
 consisting of $x_{i,j}$ and $\pd_{i,j}$ with $i,j \in I$, 
respectively. 
 Then
 $\{ \varGamma_k\}_{k=0,1,\dots,\floor{n/2} }$ 
 forms a generating system for $\PD(\Alt_n)^{\GL{n}}$,
\end{thm*}

The aim of this paper is 
to find a generating function for
$\{ \varGamma_k \}$.

Following \cite{Itoh01,KW02},
let us consider an alternating matrix 
with entries in $\PD(\Alt_n)$ given by:
\begin{equation*}
\tild{\bs{\Phi}}(u)
:=\left[
\begin{array}{cccc|cccc} 
   0     & x_{1,2} & \cdots & x_{1,n} &     &      &   &   u              
    \\
-x_{1,2} & 0       & \ddots &  \vdots  &    &      & u &           
    \\  
\vdots   &\ddots   &  0     &x_{n-1,n} &    &\adots&   &
    \\[3pt]
-x_{1,n} &\cdots   &-x_{n-1,n}& 0      & u  &      &   &
    \\[3pt]  
  \hline
       &      &      & -u   & 0      &\pd_{n-1,n}& \cdots & \pd_{1,n}
    \\
       &      &\adots&      &-\pd_{n-1,n}& 0  &\ddots  &\vdots  
    \\
       &  -u  &      &      &\vdots  &\ddots  & 0  & \pd_{1,2} 
    \\[2pt]
  -u   &      &      &      &-\pd_{1,n}&\cdots & -\pd_{1,2}& 0
\end{array}    
\right]
\end{equation*}
with $u \in \C$ a parameter.
We remark that the matrix $\tild{\bs{\Phi}}(u)$ 
(or rather, $\bs{\Phi}(u)$ given below)
naturally appears
if we regard $\GL{n}$ and $\Alt_n$ as 
a subgroup of the complex special orthogonal group $\SO_{2n}$ 
by the map \eqref{e:embed_GL_into_O} below
and the holomorphic tangent space at the origin 
of the corresponding Hermitian symmetric space of noncompact type,
respectively
(see \cite{KinvDO} for details, 
though we only deal with 
its commutative counterpart therein, 
i.e.~the principal symbol).

Our main result of this paper is the following.
Pfaffian $\Pf{\tild{\bs{\Phi}}(u)}$ of $\tild{\bs{\Phi}}(u)$
(see the next paragraph for the definition of Pfaffian)
provides a generating function for 
$\{ \varGamma_k \}$,
with coefficient being monic polynomial in $u$ of degree $n-2k$,
which is essentially equal to the Hermite polynomial,
i.e.
\begin{equation}
\label{e:main_result}
 \Pf{\tild{\bs{\Phi}}(u)}
  = \sum_{k=0}^{\floor{n/2}} \left( -\tfrac{i}{2} \right)^{n-2k}
       \hspace{-3pt} H_{n-2k}(i u) \; \varGamma_k,
\end{equation}
where $i=\sqrt{-1}$ 
and $H_m(x)$ denotes the Hermite polynomial of degree $m$.
Note that 
the minor summation formula of Pfaffian 
with commutative entries (cf.~\cite{ncpfaff08}) 
immediately implies that
the principal symbol $\sigma(\Pf{ \tild{\bs{\Phi}}(u) })$ 
of $\Pf{\tild{\bs{\Phi}}(u)}$ can be expanded as
\[
 \sigma( \Pf{\tild{\bs{\Phi}}(u) })
  = \sum_{k=0}^{\floor{n/2}} u^{n-2k} \gamma_k,
\]
where $\gamma_k$ denotes the principal symbol of $\varGamma_k$.


%
In general,
for an associative algebra $\A$
over a field $\K$ of characteristic $0$, 
which is not necessarily commutative,
Pfaffian $\Pf{\bs{A}}$ of an alternating matrix 
$\bs{A}=(A_{i,j})$, $A_{j,i}=-A_{i,j} \in \A$, 
is defined by
\allowdisplaybreaks{
\begin{align*}
 \Pf{\bs{A}}
  &= \frac{1}{2^n n!} \sum_{\sigma \in \mathfrak{S}_{2n}}
     \sgn{\sigma}\, A_{\sigma(1),\sigma(2)} A_{\sigma(3),\sigma(4)}
                     \cdots  A_{\sigma(2n-1),\sigma(2n)}
       \\
  &= \frac{1}{n!} \hspace{-10pt}
       \sum_{\begin{subarray}{c} 
          \sigma \in \mathfrak{S}_{2n} \\ 
          \sigma(2i-1) < \sigma(2i)
         \end{subarray}
       }  \hspace{-10pt}
     \sgn{\sigma}\, A_{\sigma(1),\sigma(2)} A_{\sigma(3),\sigma(4)}
                     \cdots  A_{\sigma(2n-1),\sigma(2n)}.
\end{align*}
}
(cf.~\cite{IU01}).
If the algebra $\A$ happens to be commutative,
then this reduces to:
\begin{equation*}
 \Pf{\bs{A}}
   = \sum_{\sigma}
     \sgn{\sigma}\, A_{\sigma(1),\sigma(2)} A_{\sigma(3),\sigma(4)}
                     \cdots  A_{\sigma(2n-1),\sigma(2n)},
\end{equation*}
where the summation is taken over those 
$\sigma \in \mathfrak{S}_{2n}$ satisfying
\[
 \sigma(2i-1) < \sigma(2i) \quad (i=1,2,\dots,n)
     \quad \textrm{and} \quad
 \sigma(1)<\sigma(3)<\dots<\sigma(2n-1).
\]
When dealing with Pfaffian, however,
it is sometimes convenient to consider square matrices 
alternating along the anti-diagonal,
which we call \textit{anti-alternating}
for short in this paper.
Note that a $2n \times 2n$ matrix $\bs{X}$ 
is anti-alternating 
if and only if $\bs{X} J_{2n}$ is alternating,
where $J_{2n}$ denotes the nondegenerate  
$2n \times 2n$ symmetric matrix
with $1$'s on the anti-diagonal and $0$'s elsewhere.
We simply denote $\Pf{\bs{X}J_{2n}}$ by $\Pf{\bs{X}}$
when there is no danger of confusion.
Moreover,
adopting the convention 
that $-i$ means $2n+1-i$ for $i=1,\dots,2n$,
a square matrix $\bs{X}=(X_{i,j})$ is anti-alternating
if and only if $X_{i,j}=-X_{-j,-i}$ for all $i,j$.
Thus, we will consider the anti-alternating matrix
given by
\begin{equation}
\label{e:The_matrix}
\bs{\Phi}(u):= \tild{\bs{\Phi}}(u) J_{2n}
\end{equation}
and calculate its Pfaffian in what follows.

The organization of this paper is as follows.
In Section 2,
we show that $\Pf{\bs{\Phi}(u)}$ is invariant 
under the action of $\GL{n}$.
In Section 3, 
we caluculate Pfaffian $\Pf{\bs{\Phi}(u)}$ 
and show that it provides a generating function 
for $\{ \varGamma_k \}$ 
with coefficient essentially equal to
the Hermite polynomial.

\section{Invariant differential operators}

As in the Introduction,
let $\pi$ denote the representation of $\GL{n}$
on ${\C}[\Alt_n]$ defined by
\eqref{e:GL_rep_on_polynom},
and let $M_{i,j}$ and $D_{i,j}$ denote
the multiplication operator by $x_{i,j}$ 
and the derivation  $\pd_{i,j}$, 
respectively.
The conjugation by $\pi(g)$ of them are given by the following.

%
%
\begin{lem}
\label{l:transformation_of_D_and_M}
We have
\begin{align}
 \pi(g) D_{i,j} \pi(g)^{-1}
  &= \sum_{a<b} \det(g^{a,b}_{i,j}) D_{a,b},
      \\
 \pi(g) M_{i,j} \pi(g)^{-1}
  &= \sum_{a<b} \det((g^{-1})^{i,j}_{a,b}) M_{a,b}
\end{align}
for all $g=(g_{a,b})_{a,b} \in \GL{n}$,
where $g^{a,b}_{i,j}$ denotes a $2 \times 2$ submatrix of $g$ 
whose row- and column indices are in 
$\{ a,b \}$ and $\{ i,j \}$, respectively.
\end{lem}
\begin{proof}
First, we note that 
\begin{equation}
\label{e:basic} 
g (E_{i,j}-E_{j,i}) \tp{g} 
= \sum_{a<b} \det( g^{a,b}_{i,j}) (E_{a,b}-E_{b,a}).
\end{equation}
Therefore,
setting $x=\sum_{a<b} x_{a,b} (E_{a,b}-E_{b,a})$,
we have
\begin{align*}
\pi(g) D_{i,j} \pi(g)^{-1} f(x)
&= \left. \frac{d}{d \epsilon} \right|_{\epsilon=0} 
            f(x + \epsilon g (E_{i,j}-E_{j,i}) \tp{g})
     \\
&= \left. \frac{d}{d \epsilon} \right|_{\epsilon=0} 
      f \, \biggl( \sum_{a<b} 
                  (x_{a,b} + \epsilon \det( g^{a,b}_{i,j}) )E_{a,b} 
           \biggr)
     \\
&=\sum_{a<b} \det( g^{a,b}_{i,j}) D_{a,b} f(x),
\end{align*}
and hence obtain the first formula.

As for the multiplication operator $M_{i,j}$, 
it follows from \eqref{e:basic} that
the $(i,j)$-th entry of $g^{-1} x (\tp{g})^{-1}$
equals
$\sum_{a<b} \det( (g^{-1})^{i,j}_{a,b}) x_{a,b}$.
Therefore,
\begin{align*}
\pi(g) M_{i,j} \pi(g)^{-1} f(x)
 &= \left( M_{i,j} \pi(g)^{-1} f \right)(g^{-1} x (\tp{g})^{-1} )
     \\
 &= \sum_{a<b} \det( (g^{-1})^{i.j}_{a,b}) x_{a,b}
          \left( \pi(g)^{-1} f \right)( g^{-1} x (\tp{g})^{-1} )
    \\
 &= \sum_{a<b}  \det( (g^{-1})^{i.j}_{a,b}) M_{a,b} f(x).
\end{align*}
This completes the proof.
\end{proof}

Henceforth,
we will use $x_{i,j}$ and $\pd_{i,j}$ 
to denote $M_{i,j}$ and $D_{i,j}$ 
for simplicity.


For $g \in \GL{n}$ and 
$\bs{X}=(X_{i,j}) \in \Mat{2n}(\C) \otimes \PD(\Alt_n)$,
let us denote by $\Ad_{\pi(g)}(\bs{X})$
the $2n \times 2n$ matrix 
whose $(i,j)$-th entry is given by
$\pi(g) X_{i,j} \pi(g)^{-1}$ 
for $i,j=1,\dots,2n$,
following \cite{IU01}.
Furthermore,
let 
$\SO_{2n}:=\{ g \in \GL{2n}; \tp{g} J_{2n} g=J_{2n}, \det g=1 \}$,
and $\iota $ the embedding of $\GL{n}$ into $\SO_{2n}$ 
given by
\begin{equation}
\label{e:embed_GL_into_O}
\iota: g \mapsto 
   \begin{bmatrix}
    g & 0 \\ 0 & J_n (\tp{g})^{-1} J_n
   \end{bmatrix}.
 \end{equation}

%
%
\begin{prop}
\label{p:adoint_of_Phi}
Let $\bs{\Phi}(u)$ be the matrix given by \eqref{e:The_matrix}.
Then we have
\begin{equation}
\label{e:Adjoint_of_Phi}
\Ad_{\pi(g)}(\bs{\Phi}(u))=\iota(\tp{g}) \bs{\Phi}(u) \iota(\tp{g})^{-1}
\end{equation} 
for all $g \in \GL{n}$.
\end{prop}
\begin{proof}
If we denote 
the $n \times n$ matrices 
$(\pi(g) \pd_{i,j} \pi(g)^{-1})_{i,j}$ and
$(\pi(g) x_{i,j} \pi(g)^{-1})_{i,j}$
by $\tilde{D}$ and $\tilde{M}$, respectively,
then the left-hand side of \eqref{e:Adjoint_of_Phi}
can be written as
\[
 \Ad_{\pi(g)}(\bs{\Phi}(u))
 =\begin{bmatrix}
   u 1_n & \tilde{D} J_n \\ - J_n \tilde{M} & -u 1_n
  \end{bmatrix}.
\]

On the other hand,
since
the upper-right block and the lower-left block of $\bs{\Phi}(u)$ 
can be written as $D J_n$ and $-J_n M$, respectively,
the right-hand side of \eqref{e:Adjoint_of_Phi} equals
\begin{align*} 
  &
\begin{bmatrix}  
 \tp{g} &  \\  & J_n g^{-1} J_n 
\end{bmatrix}
\begin{bmatrix}  
 u 1_n & D J_n \\ -J_n M & -u 1_n
\end{bmatrix} 
\begin{bmatrix}  
 (\tp{g})^{-1} &  \\  & J_n g J_n 
\end{bmatrix}
       \\   
  =&
\begin{bmatrix} 
 u 1_n                       & \tp{g} D g J_n  \\  
- J_n g^{-1} M (\tp{g})^{-1} & -u 1_n 
\end{bmatrix}.
\end{align*}
Now, it follows from \eqref{e:basic} that
\begin{align*}
  \tp{g} D g 
  &= \sum_{i<j} \pd_{i,j} \tp{g} (E_{i,j}-E_{j,i}) g
     \\
  &= \sum_{i<j, a<b} 
      \det( g^{a,b}_{i,j}) \pd_{a,b} (E_{i,j}-E_{j,i}),
\end{align*}
which equals the matrix $\tilde{D}$
by Lemma \ref{l:transformation_of_D_and_M}.
The same calculation shows that
\( g^{-1} M (\tp{g})^{-1} = \tilde{M} \).
Thus we obtain the proposition.
\end{proof}

As in the commutative case,
the noncommutative Pfaffian transforms 
under the action of $\GL{2n}(\K)$ 
as follows (see \cite{IU01}).
 
%
%
\begin{lem}
\label{l:adjoint_of_Pf}
Let $\bs{X}$ be an anti-alternating matrix 
with coefficient in $\A$.
For $g \in \GL{2n}(\K)$,
we have
\[
 \Pf{g \bs{X} \ttp{g}} = \det g \Pf{\bs{X}},
\]
where we set $\ttp{g}:=J_{2n} \tp{g} J_{2n}$ for brevity.
\end{lem}

By Proposition \ref{p:adoint_of_Phi} 
and Lemma \ref{l:adjoint_of_Pf}, 
we obtain the following.
%
%
\begin{cor}
\label{c:invariance_of_Pf(Phi)}
The Pfaffian $\Pf{\bs{\Phi}(u)} \in \PD(\Alt_n)$ 
is invariant under the action of $\GL{n}$.
Namely, we have 
\[
 \pi(g) \Pf{\bs{\Phi}(u)} \pi(g)^{-1} =\Pf{\bs{\Phi}(u)})
\]
for all $g \in \GL{n}$.
\end{cor}

\section{Generating function}

In this section,
we show that Pfaffian $\Pf{\bs{\Phi}(u)}$
of the matrix $\bs{\Phi}(u)$ given by \eqref{e:The_matrix}
provides a generating function for
the invariant differential operators
$\{ \varGamma_k \}$
with coefficients the Hermite polynomials,
which, combined with Corollary \ref{c:invariance_of_Pf(Phi)},
implies that 
each $\varGamma_k$ is $\GL{n}$-invariant.

As is well known,
Pfaffian is closely connected with the exterior algebra. 
Denoting by $[\pm n]$
the set $\{1,2,\dots,n, -n,\dots,-2,-1\}$,
let $V$ be a $2n$-dimensional vector space over $\K$ 
with a basis {$\{ e_i \}_{i \in [\pm n]}$}
and $\bigwedge^{\bullet} V$ the exterior algebra associated to $V$.
For $\omega,\theta \in \bigwedge^{\bullet} V$,
write the exterior product $\omega \wedge \theta$
as $\omega \theta$ for short.
Furthermore, 
let $\bigwedge^{\bullet} V \otimes \A$ be
the exterior algebra with coefficient in $\A$,
whose product is determined by 
\[
(\omega \otimes X)(\theta \otimes Y)
:=\omega \theta  \otimes XY
\]
for $\omega, \theta \in \bigwedge^{\bullet} V$ and $X,Y \in \A$
.

To an anti-alternating matrix 
$\bs{X}=(X_{i,j})_{i,j \in [\pm n]}$
with $X_{i,j} \in \A$, 
we associate 
a $2$-form $\Xi_{\bs{X}}$ defined by
\begin{equation}
\label{e:altmat_and_2form}
 \Xi_{\bs{X}} := \sum_{i,j \in [\pm n]} e_{i} e_{-j} \otimes X_{i,j}
     \in {\bigwedge\nolimits}^2 V \otimes \A.
 \end{equation}
Then the Pfaffian $\Pf{\bs{X}}$ is the coefficient of 
the volume form $e_{1} e_{2} \cdots e_{n} e_{-n} \cdots e_{-1}$
in $\Xi_{\bs{X}}^{n}$ divided by $2^n n!$:
\begin{equation}
\label{e:Pf_and_2form}
 \Xi_{\bs{X}}^n 
  = 2^n n! e_{1} e_{2} \cdots e_{n} e_{-n} \cdots e_{-1} 
     \otimes \Pf{\bs{X}}.
\end{equation}

Henceforth, 
to keep formulas concise,
for a subset $I=\{ i_1< i_2< \dots < i_k \} \subset [n]$,
put $-I:=\{ -i_k < \dots < -i_2< -i_1 \}$
and write $e_I$ and $e_{-I}$ instead of 
$e_{i_1} e_{i_2} \dots e_{i_k}$ 
and $e_{-i_k} \dots e_{-i_2} e_{-i_1}$,
respectively;
for $\omega \in \bigwedge^{\bullet} V$ and $X \in \A$,
write $\omega X$ instead of $\omega \otimes X$.

Now take $\A$ to be $\PD(\Alt_n)$,
and define 2-forms 
$\tau,\Theta_{-},\Theta_{+} \in \bigwedge^2 V \otimes \PD(\Alt_n)$ 
by
\begin{equation}
 \tau:=\sum_{i,j \in [n]} e_{i} e_{-i},
     \quad
 \Theta_{-}:=\sum_{i,j \in [n]} e_{i} e_{j} x_{i,j},
     \quad
 \Theta_{+}:=\sum_{i,j \in [n]} e_{-j} e_{-i} \pd_{i,j}.
\end{equation}
Then
$\Omega:=\Theta_{-}+2 u \tau+\Theta_{+}$
is the 2-form corresponding to $\bs{\Phi}(u)$
under $\eqref{e:altmat_and_2form}$,
and $\Pf{\bs{\Phi}(u)}$ is the coefficient of 
volume form $e_{[n]} e_{-[n]}$
in $\Omega^n$ divided by $2^n n!$.

%
%
\begin{lem}
\label{l:CR}
We have the following commutation relations among
$\tau, \Theta_{-}$ and $\Theta_{+}$:
\begin{equation}
\label{e:CR}
[\tau, \Theta_{-}]=[\tau, \Theta_{+}]=0,
   \quad
[\Theta_{+}, \Theta_{-}]=2 \tau^2.
\end{equation}
\end{lem}
\begin{proof}
These follow from easy calculation. 
For example,
we see that
\begin{align*}
[\Theta_{+},\Theta_{-}] 
 &= 4 \sum_{i<j, k<l} 
     \left(
         e_{-j}e_{-i} e_{k}e_{l} \pd_{i,j} x_{k,l}
       - e_{k}e_{l} e_{-j}e_{-i} x_{k,l} \pd_{i,j}
     \right)
     \\
 &= 4 \sum_{i<j,k<l} e_{-j}e_{-i} e_{k}e_{l} [\pd_{i,j}, x_{k,l}]
  = 4 \sum_{i<j} e_{i}e_{j} e_{-j}e_{-i},
\end{align*} 
while,
\begin{equation*}
 \tau^2
 = \sum_{i,j} e_{i}e_{-i} e_{j}e_{-j} 
 = \biggl( \sum_{i<j}+ \sum_{i>j} \biggr) \, e_{i}e_{-i} e_{j}e_{-j}
 = 2 \sum_{i<j} e_{i}e_{j} e_{-j}e_{-i}.
\end{equation*} 

\end{proof}

%
%
Let $\normalordering{\cdot}$ be 
the normal ordering in $\PD(\Alt_n)$,
i.e.~the linear map of $\Alt_n$ into itself
determined by
\[
 \normalordering{\pd_{i,j} P}
  =\normalordering{P \pd_{i,j}}
   =P \, \pd_{i,j},
     \quad
 \normalordering{x_{i,j}P} 
  =\normalordering{P x_{i,j}}
    = x_{i,j} \, P,
     \quad
 \normalordering{1}=1
\]
for any $P \in \PD(\Alt_n)$ and $i,j$. 
We extend it to 
$\bigwedge^{\bullet} V \otimes \PD(\Alt_n)$ 
canonically.
Then by definition, we obtain that
\begin{equation}
\label{e:norm_prod}
 \normalordering{(\Theta_{-}+\Theta_{+})^m}
  =\sum_{k=0}^{m} \binom{m}{k} \Theta_{-}^k \Theta_{+}^{m-k}
\end{equation}
for all $m \in \N$.

%
%
\begin{prop}
\label{p:power2normal_ordering}
Let $m$ be a nonnegative integer.
Then we have
\begin{equation}
\label{e:expansion}
  (\Theta_{-}+\Theta_{+})^m
= \sum_{k=0}^{\floor{m/2}} 
   c_{k}(m) (2\tau^2)^k \,
     \normalordering{(\Theta_{-}+\Theta_{+})^{m-2k}},
\end{equation}
where $c_k(m)$ are given by
\begin{equation}
\label{e:coeff_in_the_expansion}
 c_{k}(m)=\frac{m!}{2^k k! \, (m-2k)!}  
\end{equation}
for $k=0,1,2,\dots,\floor{m/2}$,
and $c_k(m)=0$ for $k<0$ and $k>\floor{m/2}$.
\end{prop}

We need the following lemma to prove the proposition,
though we will only use the case where $a=1$.

%
%
\begin{lem}
For nonnegative integers $a$ and $b$,
we have
\begin{equation}
\label{e:commutation_Theta+_and_Theta-}
  \Theta_{+}^a \Theta_{-}^b
 =\sum_{k=0}^{\min(a,b)}
        \frac{\dfac{a}{k} \, \dfac{b}{k}}{k!} 
          (2\tau^2)^k \Theta_{-}^{b-k}\Theta_{+}^{a-k},
\end{equation}
where,
for $z \in \C$ and $k \in \N$,
$\dfac{z}{k}$ denotes the descending factorial
$z(z-1) \cdots (z-k+1)$.
Note that $\dfac{z}{k}=0$ 
if $z \in \N$ and $k>z$.
\end{lem}
\begin{proof}
In view of the convention 
about the descending factorial,
we can assume that $a \leqsl b$ in
\eqref{e:commutation_Theta+_and_Theta-}.
Now we use induction on $a$. 
It is trivial if $a=0$.
Suppose it is true for some $a \geqsl 0$.
Then applying Lemma \ref{l:CR},
we obtain that
{\allowdisplaybreaks
\begin{align*}
\Theta_{+}^{a+1} \Theta_{-}^b
&=\sum_{k=0}^{a} \binom{a}{k} \dfac{b}{k} (2\tau^2)^k
           \Theta_{+} \Theta_{-}^{b-k} \Theta_{+}^{a-k}         
   \\
&=\sum_{k=0}^{a} \binom{a}{k} \dfac{b}{k} (2\tau^2)^k
      \left( 
           \Theta_{-}^{b-k} \Theta_{+} + [\Theta_{+},\Theta_{-}^{b-k}] 
      \right) \Theta_{+}^{a-k}
   \\
&=\sum_{k=0}^{a} \binom{a}{k} \dfac{b}{k} (2\tau^2)^k
      \left(
           \Theta_{-}^{b-k} \Theta_{+}^{a+1-k}
           + (b-k) 2\tau^2 \Theta_{+}^{b-1-k} \Theta_{-}^{a-k}
      \right)   
   \\
&=\sum_{k=0}^{a} \binom{a}{k} \dfac{b}{k} (2\tau^2)^k
      \Theta_{-}^{b-k} \Theta_{+}^{a+1-k}
  +\sum_{k=0}^{a} \binom{a}{k} \dfac{b}{k+1} (2\tau^2)^{k+1}
      \Theta_{-}^{b-1-k} \Theta_{+}^{a-k}         
   \\
&=\sum_{k=0}^{a+1} 
    \binom{a+1}{k} \dfac{b}{k} (2\tau^2)^k 
        \Theta_{-}^{b-k} \Theta_{+}^{a+1-k}.
\end{align*}
}
This completes the proof.
\end{proof}

\begin{proof}%
[Proof of Proposition \ref{p:power2normal_ordering}
]
Use induction on $m$.
There is nothing to prove when $m=0$.
Suppose that \eqref{e:expansion} is true for some $m \geqsl 0$.
Multiplying \eqref{e:expansion} by $\Theta_{-}+\Theta_{+}$ 
from the left, we obtain that
{\allowdisplaybreaks
\begin{align*}
 & (\Theta_{-}+\Theta_{+})^{m+1}
        \\
=& \sum_{k=0}^{\floor{m/2}} c_{k}(m) (2\tau^2)^k
     \sum_{s=0}^{m-2k} \binom{m-2k}{s} 
     \left( 
        \Theta_{-}^{s+1} \Theta_{+}^{m-2k-s}
       +\Theta_{+} \Theta_{-}^s \Theta_{+}^{m-2k-s}
     \right)
         \\
=& \sum_{k=0}^{\floor{m/2}} c_{k}(m)(2\tau^2)^k
    \sum_{s=0}^{m-2k} \binom{m-2k}{s}
     \left( 
        \Theta_{-}^{s+1} \Theta_{+}^{m-2k-s} 
       +\Theta_{-}^{s} \Theta_{+}^{m+1-2k-s} 
       +s 2\tau^2 \Theta_{-}^{s-1} \Theta_{+}^{m-2k-s}
     \right).
\end{align*}
Now, in the inner summation,
since $\binom{m-2k}{s-1}+\binom{m-2k}{s}=\binom{m+1-2k}{s}$,
the first and second sums equal
\begin{align*}
 & \sum_{s=0}^{m-2k} \binom{m-2k}{s}
     \left( 
          \Theta_{-}^{s+1} \Theta_{+}^{m-2k-s} 
         +\Theta_{-}^{s} \Theta_{+}^{m+1-2k-s} 
     \right)
     \\
=& \sum_{s=0}^{m+1-2k} \binom{m+1-2k}{s} 
     \Theta_{-}^{s} \Theta_{+}^{m+1-2k-s}
     \\
=& \normalordering{ (\Theta_{-}+\Theta_{+})^{m+1-2k} },
\end{align*}
while the last equals
\begin{align*}
 & \sum_{s=0}^{m-2k} \binom{m-2k}{s} s
       2\tau^2 \Theta_{-}^{s-1} \Theta_{+}^{m-2k-s} 
     \\
=& (m-2k) 2\tau^2 
     \sum_{s=0}^{m-1-2k} \binom{m-1-2k}{s}
              \Theta_{-}^{s} \Theta_{+}^{m-1-2k-s} 
     \\
=& (m-2k) 2\tau^2
     \normalordering{(\Theta_{-}+\Theta_{+})^{m-1-2k}}.
\end{align*}
Thus
\begin{align*}
 (\Theta_{-}+\Theta_{+})^{m+1}
 =& \sum_{k=0}^{\floor{m/2}} c_{k}(m) (2\tau^2)^{k}
     \normalordering{(\Theta_{-}+\Theta_{+})^{m+1-2k}}
      \\
  &+
    \sum_{k=1}^{\floor{m/2}+1} (m+2-2k)c_{k-1}(m) (2\tau^2)^{k}
     \normalordering{(\Theta_{-}+\Theta_{+})^{m+1-2k}}.
\end{align*}
}
Therefore, 
it suffices to show that
\begin{equation}
\label{e:recursive_rel}
 c_{k}(m+1)=c_{k}(m)+(m+2-2k)c_{k-1}(m),
\end{equation} 
which follows immediately from 
the definition \eqref{e:coeff_in_the_expansion}
of $c_{k}(m)$.
In fact, 
the right-hand side of \eqref{e:recursive_rel} equals
\begin{align*}
  & \frac{m!}{2^k k! \,(m-2k)!} 
   + (m+2-2k) \frac{m!}{2^{k-1}(k-1)!\,(m-2k+2)!}
      \\
 =& \frac{(m+1)!}{2^k k!\,(m-2k)!}
 = c_{k}(m+1).
\end{align*}
Hence \eqref{e:expansion} is true for $m+1$.
\end{proof}

%
%
\begin{remark}
\label{r:remark_on_CR}
Proposition \ref{p:power2normal_ordering} holds true 
in a more general situation.
Namely, let $\mathscr{A}$ be 
a noncommutative associative algebra
over an arbitrary field of characteristic $0$,
and $A, B$ two elements of $\mathscr{A}$ such that
their commutator $[A,B]:=AB-BA$ commutes 
with both $A$ and $B$:
\[
 [A,[A,B]]=[B,[A,B]]=0.
\]
Then exactly the same argument as in the proposition
yields the following formula:
\[
 (A+B)^m=\sum_{k=0}^{\floor{m/2}} 
         c_{k}(m)  \left( [A,B] \right)^{2k}
          \sum_{s=0}^{m-2k} \binom{m-2k}{s} B^s A^{m-2k-s} 
\]
with $c_{k}(m)$ given by \eqref{e:coeff_in_the_expansion}.
\end{remark}

Now we are ready.

%
%
\begin{thm}
The Pfaffian $\Pf{\bs{\Phi}(u)}$ provides 
a generating function for the $\GL{n}$-invariant 
differential operators $\{ \varGamma_k \}$: 
\[
 \Pf{\bs{\Phi}(u)}
 =\sum_{k=0}^{\floor{n/2}} a_{n-2k}(u) \varGamma_k,
\]
where $a_m(u)$ are monic polynomials in $u$
given by 
\[
a_{m}(u) = \sum_{k=0}^{\floor{m/2}}
             \frac{m!}{2^{2k} (m-2k)! k!} u^{m-2k}.
\]
for $m=0,1,2,\dots$.
\end{thm}
\begin{proof}
By Lemma \ref{l:CR} and Proposition \ref{p:power2normal_ordering},
we have
{\allowdisplaybreaks
\begin{align}
\Omega^n
&=\sum_{p=0}^{n} \binom{n}{p} (2 u \tau)^{n-p} 
            (\Theta_{-} + \Theta_{+})^{p}
     \notag \\
&= \sum_{p=0}^{n} \sum_{q=0}^{\floor{p/2}}
  \frac{n!}{(n-p)! \,q! \,(p-2q)!} (2 u  \tau)^{n-p} 
   \normalordering{ (\Theta_{-}+\Theta_{+})^{p-2q}}
     \notag \\
&=\sum_{p=0}^{n} \sum_{q=0}^{\floor{p/2}} 
   \sum_{\begin{subarray}{c}
          r,s \geqsl 0 \\ r+s=p-2q
         \end{subarray}
    }
     \frac{n!}{(n-p)! \,k! \,r! \,s!} 
       (2 u \tau)^{n-p+2q} \Theta_{-}^r \Theta_{+}^{s}.
     \label{e:omega2the_n}
\end{align}
}
Using the relations
\begin{equation*}
\Theta_{-}^r = 2^r r! \sum_{I \in \binom{[n]}{2r} }
                 e_{I} \Pf{x_I}
     \quad \textrm{and} \quad
\Theta_{+}^s = 2^s s! \sum_{J \in \binom{[n]}{2s} }
                 e_{-J} \Pf{\pd_J},
\end{equation*}
we obtain 
\begin{equation}
\Omega^n
= \sum_{p=0}^{n} \sum_{q=0}^{\floor{p/2}} \sum_{r+s=p-2q}
   \frac{n!}{(n-p)!\,q!} 2^{n-2q} u^{n-p} \tau^{n-p+2q}
     \sum_{ I \in \binom{[n]}{2r}, J \in \binom{[n]}{2s} }
     e_{-I} e_{J} \Pf{x_{I}} \Pf{\pd_{J}}.
  \label{e:omega2the_n_with_Pf}
\end{equation}
With $\tau^{n-p+2q}$ in \eqref{e:omega2the_n_with_Pf} 
expanded as
\[
 \tau^{n-p+2q}=(n-p+2q)! \sum_{K \in \binom{[n]}{n-p+2q}} e_{K}e_{-K},
\]
the only terms that survive in the summation $\sum_{K,I,J}$ 
are those corresponding to $I=J=[n] \setminus K$;
in partcular, $r=s$ and $p$ is even.
Thus the sum $\sum_{K,I,J}$ is equal to
{\allowdisplaybreaks
\begin{align*}
 & \sum_{I \in \binom{[n]}{2s}}
     e_{[n] \setminus I} e_{-[n] \setminus -I} 
     e_{I} e_{-I} \Pf{x_I} \Pf{\pd_I}
      \\
=& \sum_{I \in \binom{[n]}{2s}}
     \sgn{\begin{smallmatrix} 
           [n]  \\ [n] \setminus I, \, I
          \end{smallmatrix}} 
        e_{[n]}
     \sgn{\begin{smallmatrix} 
           -[n] \\ -[n] \setminus -I, \, -I 
          \end{smallmatrix}} 
        e_{-[n]}
       \Pf{x_I} \Pf{\pd_I}
      \\
=& e_{[n]} e_{-[n]} 
    \sum_{I \in \binom{[n]}{2s}} \Pf{x_I} \Pf{\pd_I}
\end{align*}
}
since 
$%
\sgn{\begin{smallmatrix} 
           [n]  \\ [n] \setminus I, \, I
     \end{smallmatrix}} 
 =\sgn{\begin{smallmatrix} 
           -[n] \\ -[n] \setminus -Ih, \, -I 
       \end{smallmatrix}}. 
$
Letting $p=2\nu$,
we obtain that
{\allowdisplaybreaks
\begin{align}
\Pf{\bs{\Phi}(u)}
=& \sum_{\nu=0}^{\floor{n/2}} \frac{u^{n-2\nu}}{(n-2\nu)!}
   \sum_{s=0}^{\nu}
    \frac{(n-2s)!}{(\nu-s)! \, 2^{2(\nu-s)}} \varGamma_{s}
      \notag  \\
=& \sum_{s=0}^{\floor{n/2}} 
    \sum_{\nu=s}^{\floor{n/2}} 
       \frac{ (n-2s)! }{ (n-2\nu)! (\nu-s)! \, 2^{2(\nu-s)} } \,
         u^{n-2\nu} \, \varGamma_{s}
      \notag \\
=& \sum_{s=0}^{\floor{n/2}} 
       \sum_{r=0}^{\floor{n/2}-s} 
         \frac{ (n-2s)! }{(n-2s-2r)! r! \, 2^{2r} } \,
          u^{n-2s-2r} \, \varGamma_{s}
    \notag \\
=& \sum_{s=0}^{\floor{n/2}} a_{n-2s}(u) \, \varGamma_s.
    \notag
\end{align}
}
This completes the proof.
\end{proof}

The polynomials $a_m(u)$ are essentially equal to 
the Hermite polynomials given by
$H_m(x)=(-1)^m e^{x^2} \left(\frac{\d}{\d x}\right)^{m} e^{-x^2}$.
In fact, 
it is well known that the generating function for $H_m(x)$ 
is given by
\[
 e^{2tx-t^2}=\sum_{m=0}^{\infty} \frac{t^m}{m!} H_m(x),
\]
from which one can derive that
\[
 H_m(x)=m! \sum_{k=0}^{\floor{m/2}} 
            \frac{(-1)^k 2^{m-2k}}{k!(m-2k)!}x^{m-2k}.
\]
Therefore,
$a_m(u)=\left(-\tfrac{\sqrt{-1}}{2}\right)^m H_m(\sqrt{-1} u)$,
and we obtain \eqref{e:main_result}. 



\bibliographystyle{amsplain}
\bibliography{rep}

\nocite{ncpfaff08}
\nocite{GW98}

\end{document}